\documentclass[12pt, reqno]{amsart}
\usepackage{hyperref, xy}
\usepackage{amssymb}

\theoremstyle{plain}
\newtheorem*{thm*}{Theorem}
\newtheorem{thm}{Theorem}[section]

\newtheorem{lemma}[thm]{Lemma}
\newtheorem{cor}[thm]{Corollary}
\newtheorem{qn}[thm]{Question}
\theoremstyle{definition}

\theoremstyle{remark}

\newtheorem*{rem*}{Remark}
\newtheorem*{example*}{Example}
\numberwithin{equation}{section}

\begin{document}
\title[
Enveloping Algebras and Sklyanin Algebras]{ Faithful Cyclic
Modules for Enveloping Algebras and Sklyanin Algebras.}
\author{Ian M. Musson}
\address{
\newline University of Wisconsin-Milwaukee
\newline Department of Mathematical Sciences \newline P. O. Box 0413
\newline Milwaukee, WI 53211-1512}
\email{musson@uwm.edu} {\it Dedicated to Don Passman on the
occasion of his $65^{th}$ birthday.}
\thanks{Research partly supported by a grant from the NSA}
\date{\today}
\maketitle
\begin{abstract}
 Let $U$ be the enveloping algebra of a finite dimensional nonabelian Lie
 algebra $\mathfrak{g}$ over a field of characteristic zero.  We show that there
 is an open nonempty open subset $X$ of $U_1 = \mathfrak{g}\oplus K$ such that
 $U/Ux$ is faithful for all $x \in X$.  We prove similar results
 for homogenized enveloping algebras and for the three dimensional
 Sklyanin algebras at points of infinite order.  It would be
 interesting to know if there is a common generalization of these
 results.
\end{abstract}
 A prime Noetherian ring $U$ is {\it bounded} if every
essential left ideal contains a nonzero two sided ideal.  We say
that $U$ is {\it fully bounded Noetherian} (FBN) if every prime
image of $U$ is bounded.  The main examples of FBN rings are
Noetherian rings satisfying a polynomial identity (PI), see
\cite{AmSm}.
 If a filtered ring $U = \cup_{n \geq 0}U_n$ is not
bounded, it is reasonable to ask if we can find a regular element
$x \in U_1$ such that $U/Ux$ is faithful. In addition we might
hope that most elements $x \in U_1$ have this property. In this
paper we show that this is the case for enveloping algebras of
nonabelian Lie algebras in characteristic zero. We prove similar
results
 for homogenized enveloping algebras and for the three dimensional
 Sklyanin algebras.  Since these graded algebras are defined in rather different ways, we might expect that there is a more general underlying result.  We remark that the proof for nonnilpotent Lie algebras is rather similar to the Sklyanin algebra case, but for nilpotent Lie algebras a different argument is required.

 \indent  For more information and references on enveloping
algebras satisfying a PI the reader may consult Chapter 4 of
\cite{BMPZ}. Group algebras satisfying a PI are
discussed at length in the masterful book of Passman, [P]. \\
\indent  My interest in the issues raised in this paper was
stimulated  by a question of Jorg Feldvoss at the conference in
Madison.  Jorg asked whether an enveloping algebra which is FBN
must necessarily be commutative, in characteristic zero.  I thank
Jorg for his question.  In addition I thank Michaela Vancliff for
pointing me towards some relevant results in the literature.

\section{Enveloping Algebras.}

We assume that $\mathfrak{g}$ is a finite dimensional Lie algebra
over an algebraically closed field $K$ of characteristic zero. Let
$\{U_n\}_{n \geq 0} $ be the filtration on $U = U(\mathfrak{g})$
given by $U_0 = K, U_1 = \mathfrak{g}\oplus K$ and $U_n = U_1^n$
for $n \geq 2 .$  We write  $S = Gr U = \oplus_{n \geq 0}
U_n/U_{n-1}$ for the associated graded ring.
\begin{thm}\label{ea}
If $\mathfrak{g}$ is nonabelian there is a nonempty open subset
$X$ of $U_1$ such that $U/Ux$ is faithful for all $x \in X$.
\end{thm}

\begin{cor}
If $U$ is bounded then $\mathfrak{g}$ is abelian.
\end{cor}
\noindent{\bf Proof.}  This follows immediately from Theorem \ref{ea}.\\

The proof of Theorem \ref{ea} will be given as a series of lemmas.
The subset $X$ in the Theorem has the form $X = X' \oplus K,$
where $X'$ is a nonempty open subset of $\mathfrak{g}.$ First we
suppose that $\mathfrak{g}$ is not nilpotent. Then by Engel's
Theorem the set of $ad$-nilpotent elements of $\mathfrak{g}$ is a
proper closed subset, and we let $X'$ be its complement.

If $x' \in X'$, then $ad \; x'$ has a nonzero eigenvalue $\lambda$
on $\mathfrak{g}$, so there exists $y \in \mathfrak{g}$ such that
$[x',y] = \lambda y$.  Replacing $x'$ by $\lambda^{-1} x'$ we can
assume that $\lambda = 1$.  If $\mu \in K,$ and $x = x' + \mu,$
then for all $n \geq 0$ we have in $U$
\begin{equation} \label{X2}
(x-n)y^n = y^n x
\end{equation}
Using the PBW theorem we see that $y^n \notin Ux.$  If $m$ is an
element of a left $U$-module $M,$ we write $ ann_U m$ for the
annihilator of $m$ in $U.$
\begin{lemma}\label{x1}
Let $M = U/Ux$, and denote the image of $u \in U$ in $M$ by
$\overline{u}$.  Then
\[ ann_U \overline{y}^n = U(x-n). \]
\end{lemma}
\noindent{\bf Proof.} It follows from (\ref{X2}) that $U(x-n)
\subseteq ann_U \overline{y}^n$.  If $u \in U_m \cap ann_U
\overline{y}^n$, we show by induction on $m$ that $u \in U(x-n)$.
Define a filtration $\{M_p\}$ on $M$ by setting
\[ M_p = (U_p + Ux)/Ux \]
for $p \geq 0$.
Then the associated graded module $Gr M = \oplus_{p \geq 0}
M_p/M_{p-1}$ is isomorphic to $S/Sx'$, so every nonzero element of
$Gr \; M$ has annihilator $Sx'$.  Hence there exists $v \in
U_{m-1}$ such that
\[ u' = u - v(x-n) \in U_{m-1} . \]
Since $u' \overline{y}^n = 0$ the result follows by induction.
\begin{lemma} \label{x4}
We have $\cap_{n \geq 0} U(x-n) = 0.$
\end{lemma}
\noindent{\bf Proof.} Let $x_1, \ldots, x_m, x'$ be a basis for
$\mathfrak{g}$, and let $W$ be the subspace of $U$ spanned by all
elements of the form
\[ x^{a_1}_1 \ldots x^{a_m}_m \]
with $a_1, \ldots, a_m \in \mathbb{N}$.  It follows easily from
the PBW theorem that any element $u$ of $U$ can be written
uniquely in the form
\[ u = \sum^N_{i=0} w_i x^i \]
with $w_0, \ldots, w_N \in W$.  If $u \in U(x-n)$, it is easy to
see that $\sum^N_{i=0} n^i w_i = 0$.  If $u \in \cap_{n \geq 0}
U(x-n)$ we conclude from the nonvanishing of the Vandermonde
determinant that all the $w_i$ are zero.\\

\indent If $\mathfrak{g}$ is not nilpotent Theorem \ref{ea}
follows from Lemmas \ref{x1} and \ref{x4}.  Now suppose that
$\mathfrak{g}$ is nilpotent with center $\mathfrak{z}$, and let
$Z$ be the center of $U$.  Theorem \ref{ea}  in this case follows
from the next result.
\begin{lemma}\label{X3}
Assume $\mathfrak{g}$ is nilpotent and nonabelian, and set $X' = \mathfrak{g}
\backslash
 \mathfrak{z}, X = X' \oplus K.$  If $x \in X$, then $M = U/Ux$ is a faithful
 $U$-module.
\end{lemma}
\noindent {\bf Proof.}  If $I = ann_U M$ is nonzero then by
\cite{D}, Proposition 4.7.1, $I \cap Z$ is nonzero.  It suffices
to show that $Ux \cap Z = 0$ since then $M$ contains the free
$Z$-module $(Z + Ux)/Ux$.

\indent Suppose that $x'\in X, \mu \in K$ and $x = x' + \mu.$
Choose $y \in \mathfrak{g}$ such that $[y,x'] = [y,x] \neq 0$. We
extend $ady$ to a locally nilpotent derivation of $U$. Assume $u
\in U$ is nonzero, and let $i,j$ be the least integers such that
\[ (ady)^{i+1}(u) = 0 = (ady)^{j+1}(x) . \]
Then $i \geq 0$ and $j \geq 1$.  If $n = i + j$, then
\[ (ady)^n (ux) = \left( \begin{array}{c} n\\i \end{array}\right)
(ady)^i (u) (ady)^j (x), \] and this is nonzero since $U$ is a
domain.  It follows that $ux \not\in Z$.
\section{Graded Algebras.}

Now we suppose that $A$ is a $3$ dimensional Sklyanin algebra,
that is $3$-dimensional generic regular algebra of type $A$ with
$3$ generators and $3$ quadratic relations, see \cite{ATV1},
Section 4.13 or \cite{ArSc}, (10.14).  To construct $A,$ consider
a nonsingular cubic curve $E$ in $\mathbb{P}^2 = \mathbb{P}(V^*).$
We use the symbols $\oplus, \ominus$ to addition or subtraction
using the group law on $E,$ since we also need to consider
addition of divisors. We suppose the identity $0$ in $E$ is a
flex, so that $P\oplus  Q \oplus R = 0$ if and only if $P,Q$ and
$R$ are colinear.  Fix a $p$ point on $E$ and let $\sigma$ be the
an automorphism of $E$ defined by $\sigma(P) = P \ominus p.$

Given this data, we define $A$ as follows. Let $T(V) =
\bigoplus_{n \geq 0} V^{\otimes n}$ be the tensor algebra on $V.$
For ${n > 0}$ let $E_n$ be the subset of $E \times \ldots \times
E$ ($n$ copies) consisting of all $n$-tuples $(e,\sigma(e), \ldots
,\sigma^{n-1}(e))$ with $e \in E,$ and set
\[\mathcal{R}_n = \{ f \in V^{\otimes n}|f \mbox{ vanishes on} \; E_n\}. \]

For ${n \geq 0}$ let $I_n$ be the ideal of $T(V)$ generated by the
subspaces $\mathcal{R}_i$ for $2 \leq i \leq n.$ We set $A =
T(V)/I_2.$ Then $A = \bigoplus_{n \geq 0} A(n)$ where $A(n)$ is
the image of $V^{\otimes n}$ in $A.$ Let $L$ be the line bundle
associated to the embedding of $E$ into $\mathbb{P}^2.$ In
\cite{ATV1} and \cite{ATV2} the algebra $A$ is denoted by
$A(\mathcal{T})$ where  $\mathcal{T}$ is the triple
$(E,\sigma,L).$
 By \cite{ATV1}  Theorem 6.8, there is an element
$g \in A(3)$ such that $gA = Ag$ and $A/Ag$ is isomorphic to the
twisted homogeneous coordinate ring $B,$  (denoted
$B(\mathcal{T})$ in \cite{ATV1}).  It follows that $I_n = I_3$ for
all ${n \geq 3}.$ The algebra $A$ is a domain by \cite{ATV2}
Section 3, and is finite module over its center if and only if $p$
is a point of finite order, \cite{ATV2} Theorem 7.1. Hence by
Posner's theorem, \cite{MR} Theorem 13.6.5, $A$ is a PI algebra if
and only if $p$ is a point of finite order. It is well known that
$B$ is a domain. To see this we can adapt the proof of \cite{V},
Lemma 3.5.

\begin{thm} \label{Y11} If $p$ is a point of infinite order on $E$
then $A/AL$ is a faithful $A$-module for all $L \in A(1).$
\end{thm}

Again the proof will be given as a series of lemmas. Equation
(\ref{LP1}) and Lemmas \ref{Y2} and \ref{Y12} are analogs of
equation (\ref{X2}) and Lemmas \ref{x1} and \ref{x4} respectively.

 We refer to nonzero elements of $A(1) = V$ as {\it
lines}.
A {\it divisor} on $E$ is a finite formal sum of the form $\sum_{P
\in E} n_pP$ with $n_p \in \mathbb{Z}.$  If $D = \sum_{P \in E}
n_pP$ we set $\sigma_*(D) = \sum_{P \in E} n_p\sigma(P).$ The
divisor $D$ is {\it effective} if $n_p \geq 0$ for all $P \in E.$
If $D$ is effective we say that $P$ {\it occurs in} $E$ if $n_p >
0.$

If $F \in V^{\otimes n}$ does not vanish on $E_n$ we write
$Div_n(F)$ for the divisor on $E$ cut out by $F,$ that is
\[Div_n(F) = \sum \{e \in E |
F(e,\sigma(e), \ldots ,\sigma^{n-1}(e)) = 0 \}\] where zeroes of
$F$ are counted with multiplicities. If $F \in V^{\otimes n}$ then
$Div_n(F)$ depends only on the image of $F$ modulo
$\mathcal{R}_n.$  Hence if $f \in A(n),$ and $F \in V^{\otimes n}$
is a preimage of $f$ in $V^{\otimes n},$ we set $Div_n(f) =
Div_n(F).$
 For $F \in V$ we write $Div(F)$ in place of $Div_1(F).$

Note that if $F \in V^{\otimes (n-1)},$ and $L \in V,$ then
\[Div_n(F\otimes L) = Div_{n-1}(F) + \sigma_*^{1-n}(Div(L))\] and
\[Div_n(L \otimes F) =
Div(L) + \sigma_*^{-1}(Div_{n-1}(F)).\]
Let $L_0 = L$ and suppose
that $Div(L_0) = P + Q + R.$ Let $\mathbb{N}^* = \mathbb{N}
\backslash \{0 \} $ be the set of positive integers.  Since $P, Q,
R$ are colinear, and $p$ has infinite order, we cannot have $3P,
3Q, 3R \in 3\mathbb{N}^*p,$ and we may assume that $ 3R \notin
3\mathbb{N}^*p.$ For all $n
> 0$ fix a line $L_n$ through $P\oplus np$ and $Q\oplus np$
and $R \ominus 2np.$

Choose a line  $M_0$ through $R$ such that $M_0$ is not a scalar
multiple of $L_0.$ Then $Div(M_0) = R + S + T,$ where $\{P, Q\}
\cap \{S, T\} = \emptyset.$

\begin{lemma} \label{Y1} For all $n > 0$
 there are lines $M_n$ such that

\begin{itemize}
\item[{ (1)}] $Div(M_n) =  (R \ominus 2np) + (S\oplus np) +
(T\oplus np).$
 \item[{ (2)}]
 $L_{n}M_{n-1} = M_{n}L_{n-1}$
holds in $A(2).$
\end{itemize}
\end{lemma}
\noindent{\bf Proof.}  Since the points $R \ominus 2np,\; S\oplus
np$ and $T\oplus np$ sum to zero using the group law on $E,$ there
are lines $M_n$ satisfying (1). It suffices to show that (2) holds
when  $n = 1.$
Since $Div_2(L_{1} \otimes M_0)$ and $Div_2(M_{1} \otimes L_0)$
are both equal to
\[(P \oplus p) + (Q\oplus p) + (R \ominus 2p) + (R \oplus p)  +
(S\oplus p) + (T\oplus p),\] it follows that $L_1 {\otimes } M_0 -
\lambda M_1 {\otimes } L_0 \in \mathcal{R}_2$ for some nonzero
scalar $\lambda.$  Replacing $M_1$ by $\lambda^{-1}M_1$ we obtain
the result.\\

Now set $N_1 = M_0,N'_1 = M_1,$ and for $i > 1$ define $N_i,N'_i$
inductively by $N_i = M_{i-1}N_{i-1}$ and $N'_i = M_iN'_{i-1}.$ By
Lemma \ref{Y1} and induction, we have
 for $n \geq 1$  \begin{equation}\label{LP1}
L_nN_n = N'_n L.
\end{equation}

We write $GK(M)$ for the Gelfand-Kirillov dimension of the graded
$A$-module $M.$ This can be defined as the order of the pole of
the Hilbert series of $M,$ see \cite{ATV2} Section 2. We say that
$M$ is {\it critical} if $GK(N) < GK(M)$ for every proper graded
factor module $N$ of $M.$ Since the Hilbert series is additive on
short exact sequences we see that if $M$ is  critical then any
nonzero graded submodule of $M$ has GK-dimension equal to $GK(M).$

\begin{lemma} \label{Y2}
If $\overline{N}_n$ is the image of $N_n$ in $A/AL,$ then  $ann_A
\overline{N}_n = AL_n.$
\end{lemma}
\noindent{\bf Proof.} We show first that $N_n \notin AL.$  From
Lemma \ref{Y1},  it follows that
\[ Div(N_n) =  n(S \oplus (n-1)p) +
n(T\oplus (n-1)p) + \sum^{n-1}_{i=0}
R \oplus (3i -2(n-1))p.\]

On the other hand if $F \in A(n-1),$ then
\[ Div_n(FL) = Div_{n-1}(F) +  (P \oplus (n-1)p) +
(Q \oplus (n-1)p) + (R \oplus (n-1))p.\] If $N_n = FL,$ then since
$\{P, Q\} \cap \{S, T\} = \emptyset,$ it would follow that $P  = R
\ominus 3ip,$ and $Q  = R \ominus 3jp,$ where $0 < i, j \leq n-1.$
Then since $P,Q$ and $R$ are colinear we would have $ 3R = 3(i +
j)p \in 3\mathbb{N}^*p,$ a contradiction.\\
\indent By equation (\ref{LP1}),   $AL_n \subseteq ann_A
\overline{N}_n.$ It follows from \cite{ATV2}, Proposition 6.1,
that $A/AL$ is a critical $A$-module with GK-dimension 2.
Therefore the submodule $A \overline{N}_n$ of $A/AL$ has
GK-dimension 2, and any proper factor module of  $A/AL_n$ has
GK-dimension less than 2. Hence the natural map from $A/AL_n$ onto
$A \overline{N}_n$ must be an isomorphism, and this implies the
result.

\begin{lemma} \label{Y12}
$\cap_{n \geq 0}AL_n = 0.$
\end{lemma}
\noindent{\bf Proof.} We first show that $\cap_{n \geq 0}AL_n
\subseteq (g).$  Suppose that $f \in \cap_{n \geq 0}AL_n.$ We can
assume that $f \in A(m)$ for some $m.$ Let $F \in V^{\otimes m}$
be a preimage of $f$ in  $T(V).$ If $f \notin (g),$  then $F
\notin \mathcal{R}_m,$ so $F$ does not vanish on $E_m$ and we can
consider the effective divisor $Div_m(F).$ For all $n$ we can
write $f = G_nL_n$ for some $G_n \in A(m-1)$ but then $Div_m(f)$
contains $\sigma^{1-m}(P \oplus np) = P \oplus (n + m -1)p$ for
all ${n \geq 0}.$ However the divisor $Div_m(f)$ cannot contain
more than $3m$ points.  Since $p$ is a point of infinite order
this is a contradiction.

Next we show  that $\cap_{n \geq 0}AL_n \subseteq (g^i)$ for all
$i \geq 1$. Any nonzero element of $(g^i)$ has degree at least
$3i,$ so $\cap_{i \geq 1}(g^i) = 0$ and this will finish the
proof. Suppose that $x = a_nL_n$ with $a_n \in A,$ for all $n \geq
0.$ By induction we assume that $a_n \in (g^i)$ for all $n$.  Then
$x = g^ib_nL_n,$ for some $b_n \in A,$ and $y = b_nL_n$ is
independent of  $n$ since $A$ is a domain. By the first part of
the proof $y \in (g).$ Since $B = A/(g)$ is a domain and $L_n
\notin (g)$ by comparing degrees, we see that $ b_n \in (g),$ and
then $ a_n \in (g^{i+1})$ as required.\\

Theorem \ref{Y11} follows at once from Lemmas \ref{Y2} and
\ref{Y12}.

\section{Concluding Remarks.}
The enveloping algebra result can be formulated in terms of graded
rings.  In general if $\{ U_n\}_{n \geq 0} $ is a filtration on a
ring $U$, the {\it Rees
 ring} $\widetilde{U}$ of this filtration is the
 graded subring $\widetilde{U} =
\bigoplus_{n \geq 0} \widetilde{U}(n)$ of the polynomial ring
$U[z]$ with $\widetilde{U}(n) =
 U_nz^n.$

\begin{lemma} \label{RR}
Suppose that $x \in U_n$ and set $\widetilde{x} = xz^n.$ Assume
that $Gr U$ is a domain. Then $U/Ux$ is a faithful $U$-module if
and only if $\widetilde{U}/ \widetilde{U}\widetilde{x}$ is a
faithful $\widetilde{U}$-module.
\end{lemma}
\noindent{\bf Proof.} Assume that $U/Ux$ is faithful.  Since the
annihilator of $\widetilde{U}/\widetilde{U}\widetilde{x}$ is a
graded ideal, it is enough to show that if $u = u_mz^m \in
\widetilde{U}(m)$, and $u \widetilde{U} \subseteq
\widetilde{U}\widetilde{x}$, then $u = 0$.  However
\[u\widetilde{U}(p) \subseteq \widetilde{U}\widetilde{x} \cap \widetilde{U}(p + m) =
\widetilde{U}(p + m-n)\widetilde{x},\] implies that $u_mU_p
\subseteq U_{p + m-n}x.$ In particular $u_mU \subseteq Ux,$ so $u
= 0.$
  \indent Conversely, suppose that
  $\widetilde{U}/\widetilde{U}\widetilde{x}$ is faithful.  If $a
  \in U_m \backslash U_{m-1}$ we set $\deg a = m$.  Since $Gr U$ is a domain
  \[ \deg(ab) = \deg a + \deg b. \]
 Since $x \in
 U_n$  we have $\deg x = n - i$ with $i \geq 0.$ Suppose $a \in ann_U U/Ux$ and $\deg a = m$.  If $b \in U$,
and $\deg b = p$, then
\[ ab = cx \]
for some $c \in U$ with $\deg c = m + p - n + i$.  If $u =
az^{m+i}$ then
\[ u\cdot bz^{p} = cz^{p+m-n+i} \cdot xz^n \in
\widetilde{U}\widetilde{x} . \]
 Thus, $u = 0$ by hypothesis, so $a = 0$.\\

   \indent When $U = U(\mathfrak{g})$ is an enveloping
algebra with the standard filtration $\{U_n\}_{n \geq 0},$ the
Rees ring is known as the {\it homogenized enveloping algebra} and
denoted by $H(\mathfrak{g}).$ We have $H(\mathfrak{g})(1) = Kz
\oplus \mathfrak{g}z.$

\begin{cor} 
Let $\mathfrak{g}$ be a finite dimensional nonabelian Lie algebra,
and set $H = H(\mathfrak{g}).$ Let $X$ be the open subset of $U_1$
given in Theorem \ref{ea}.  If $\widetilde{X} = Xz \subseteq H(1)$
then $H/HL$ is a faithful $H$-module for all $L \in
\widetilde{X}.$
\end{cor}
\noindent{\bf Proof.}  Combine Theorem \ref{ea} and Lemma
\ref{RR}.\\

Our results prompt the following
\begin{qn}
Suppose that $A$ is a connected graded algebra over a field and
that $A$ is generated $A(1)$ which is finite dimensional. If $A$
is not FBN, is the set
\[ \{x \in A(1)| A/Ax \; \mbox{is faithful} \} \]
Zariski dense in $A(1)?$\\
\end{qn}
To study this question the following result should be useful. We
give $Spec A$ has the Jacobson topology, so that the closed sets
have the form
\[ V(I) = \{ P \in Spec A | I \subseteq P \} \]
where $I$ is an ideal of $A$.
\begin{thm} \label{closed}
Suppose that $A = \bigoplus_{n \geq 0}A(n)$ is a
graded $K$-algebra which is a domain. Assume that  $dim_K A(n) <
\infty$ for all $n,$ and that if $L \in A(1)$ the $A$-module $M =
A/AL$ has prime annihilator.  Then the map \[ A(1) \longrightarrow
Spec A,  \; L \mapsto ann_A A/AL\] is continuous.
\end{thm}
We remark that the hypothesis that $A/AL$ has prime annihilator
for all $L \in A(1)$ holds for regular algebras of dimension $3$
by \cite{ATV2} Propositions 2.30 and 6.1.  The proof of the
theorem is based on the next lemma.

\begin{lemma} 
Keep the hypothesis of the theorem.
If $F \in A(n)$ for $n \geq 1$ then the set $X = \{ L \in A(1)|F
\in AL\}$ is Zariski closed.
\end{lemma}
\noindent
 {\bf Proof.} Let $U = A(n-1), V = A(1)$ and $W = A(n)$.  If $u
 \in U$ is nonzero, let $[u] \in \mathbb{P}(U)$ denote the line
 through $u$.  We use similar notation for elements of $V,W$.
 Since $A$ is a domain, there is a well defined morphism of
 varieties
 \[ \pi : \mathbb{P} = \mathbb{P}(U) \times \mathbb{P}(V)
 \longrightarrow \mathbb{P}(W) \]
 given by
 \[ ([u], [v]) \longrightarrow [uv]. \]
 Thus the set $\pi^{-1}([F])$ is closed in $\mathbb{P}$.  Now if $p:
 \mathbb{P} \longrightarrow \mathbb{P}(V)$ is the projection onto the
 second factor, it follows from
 \cite{H}, Theorem 3.12 that $Y = p(\pi^{-1}([F]))$ is closed in $\mathbb{P}(V)$.
It is easy to see that
 \[ X = \{ L \in A(1)|[L] \in Y\} \]
  and it follows that $X$ is closed in $A(1)$.\\
\\
{\bf Proof of Theorem \ref{closed}. } If $I$ is an ideal of $A$,
then
\begin{eqnarray*}
\pi^{-1}(V(I)) & = & \{ L \in A(1) | \pi(L) \in V(I) \} \\
                & = & \{ L \in A(1) | I \subseteq AL \} \\
                & = & \cap_{F \in I} \{ L \in A(1) | F \in AL \}
\end{eqnarray*}
and this is an intersection of closed sets by the Lemma.

\begin{qn}
For a group algebra analog of the problems considered in this
paper set
\[X = \{g
\in G| KG/KG(g - 1) \; \mbox{is faithful} \}.
\] What can be said about $X$?\\
\end{qn}


\end{document}